\input amstex
\magnification=1200
\font\cyr=wncyr10
\font\cyb=wncyb10
\font\cyi=wncyi10
\font\cyre=wncyr8
\font\cyie=wncyi8
\documentstyle{amsppt}
\NoRunningHeads
\NoBlackBoxes
\define\mYB{\operatorname{\widetilde mYB}}
\define\End{\operatorname{End}}
\define\Mat{\operatorname{Mat}}
\define\PG{\operatorname{\Pi\Gamma}}
\define\soa{\operatorname{\frak s\frak a}}
\define\sla{\operatorname{\frak s\frak l}}
\define\ad{\operatorname{ad}}
\topmatter
\title Topics in isotopic pairs and their representations. III.
Bunches of Lie algebras and modified classical Yang-Baxter equation
\endtitle
\author\bf D.V.Juriev\footnote{\ This is an English translation of the
original Russian version, which is located at the end of the article as an
appendix. In the case of any differences between English and Russian
versions caused by a translation the least has the priority as the original
one.\newline}
\endauthor
\date q-alg/9708027\enddate
\affil\eightpoint Research Center for Mathematical Physics and Informatics
``Thalassa Aitheria'',\linebreak
ul.Miklukho-Maklaya 20-180, Moscow 117437 Russia.\linebreak
E-mail: denis\@juriev.msk.ru
\endaffil
\endtopmatter
\document This short note is a brief comment to certain aspects of the theory
of classical $r$-matrix [1] and bihamiltonian formalism [2], which motivations
lie in constructions of articles [3].

\head\bf 1. Linear bunches of Lie algebras and Lie $\mYB$-algebras\endhead

This paragraph is devoted to the anlysis of some algebraic structures
associated with linear bunches of Lie algebras.

\subhead 1.1. $\Gamma$-bunches of Lie algebras and tangent bracket
\endsubhead

\definition{Definition 1A} The {\it (one-parametric) bunch of Lie algebras\/}
is a linear space $\frak g$ supplied by a family of Lie brackets
$[\cdot,\cdot]_\lambda$ ($\lambda\in\Bbb R$). The Lie algebra from $\frak g$,
which is defined by the bracket $[\cdot,\cdot]_\lambda$, is denoted by
$\frak g_\lambda$. A bunch of Lie algebras $\frak g$ is called the
{\it $\Gamma$-bunch} if there are defined homomorphisms $R_\lambda$ from
the Lie algebras $\frak g_\lambda$ into the Lie algebra $\frak g_0$,
otherwords, the equality
$$R_\lambda[X,Y]_\lambda=[R_\lambda X,R_\lambda Y],$$
where $[\cdot,\cdot]=[\cdot,\cdot]_0$, holds for all $X$ and $Y$ from
$\frak g$.
\enddefinition

Below there will be considered only smooth bunches, i.e. the bunches,
which Lie brackets form a smooth family, and smooth $\Gamma$-bunches,
for which the family of homomorphisms $R_\lambda$ is also smooth.

\proclaim{Proposition 1} Let $\frak g$ be a smooth $\Gamma$-bunch of
Lie algebras. Define the tangent bracket $[\cdot,\cdot]_R$ as
$$[X,Y]_R=\left.\frac{d[X,Y]_\lambda}{d\lambda}\right|_{\lambda=0}.$$
Then
$$[X,Y]_R=[RX,Y]+[X,RY]-R[X,Y],$$
where
$$R=\left.\frac{dR_\lambda}{d\lambda}\right|_{\lambda=0}.$$
\endproclaim

\remark{Remark 1} The bracket $[\cdot,\cdot]_R$ is trivial (is identical
zero) if and only if $R$ is a derivative of the bracket $[\cdot,\cdot]$.
\endremark

\remark{Remark 2} The bracket $[\cdot,\cdot]_R$ is a Lie one if and only if
$$\bigl([B(X,Y),Z]+B([X,Y],Z)\bigr)+c.p.=0,$$
where
$$B(X,Y)=R[RX,Y]+R[X,RY]-[RX,RY]-R^2[X,Y].$$
The expanded expression for $B(X,Y)$ may be compactified into the binomial
for\-mu\-la
$$B(X,Y)=R[X,Y]_R-[RX,RY].$$
\endremark

\subhead 1.2. Linear $\Gamma$-bunches of Lie algebras and the modified
Yang-Baxter equation\endsubhead Note that the bracket $[\cdot,\cdot]$
and the tangent bracket $[\cdot,\cdot]_R$ in any $\Gamma$-bunch
of Lie algebras always obey the identity
$$\bigl([[X,Y]_R,Z]+[[X,Y],Z]_R\bigr)+c.p.=0.$$
Hence, if the tangent bracket $[\cdot,\cdot]_R$ is a Lie one, then
it is compatible with the bracket $[\cdot,\cdot]$. This circumstance
motivates to consider the linear $\Gamma$-bunches of Lie algebras.

A bunch of Lie algebras $\frak g$ is called linear if the Lie brackets
$[\cdot,\cdot]_\lambda$ from a linear family, otherwords,
$[\cdot,\cdot]_\lambda=(1-\lambda)[\cdot,\cdot]+\lambda[\cdot,\cdot]_1$.
A $\Gamma$-bunch of of Lie algebras is called linear if not only Lie
brackets but also the homomorphisms $R_\lambda$ form a linear family,
i.e. $R_\lambda=1+\lambda R$. The condition of linearity of
the family of homomorphisms $R_\lambda$ is natural, because for any
nonlinear family $R_\lambda$ its linear part $1+\lambda R$
($R=\left.\frac{dR_\lambda} {d\lambda}\right|_{\lambda=0}$) defines the same
linear family of brackets $[\cdot,\cdot]_\lambda$.

\proclaim{Theorem 1A} For a linear $\Gamma$-bunch of Lie algebras $\frak g$
the tangent bracket $[\cdot,\cdot]_R$ (which is a Lie one in this case)
obeys the following condition
$$R[X,Y]_R=[RX,RY],$$
otherwords,
$$R[RX,Y]+R[X,RY]=[RX,RY]+R^2[X,Y].$$
\endproclaim

The theorem means that the condition that the bracket $[\cdot,\cdot]_R$ is a
Lie one from remark 2 is factorized in case of linear $\Gamma$-bunches and is
reduced to the claim of the identical equality of the expression $B(X,Y)$ to
zero. Note that the fact that the bracket $[\cdot,\cdot]_R$ is a Lie one in
view of its compatibility with $[\cdot,\cdot]$ means that the bunch
$[\cdot,\cdot]+\lambda[\cdot,\cdot]_R$ is a linear bunch of Lie algebras,
however, it is not obligatory a $\Gamma$-bunch, therefore the condition
of remark 2 is not factorized generally.

The least identity of theorem 1A coincides with the modified classical
Yang-Baxter equation [1] if $R^2=1$ that holds for the most important
examples con\-si\-dered in [1]. In this case as it was shown in [1] the
bracket $[\cdot,\cdot]'_R$ defined as
$$[X,Y]'_R=[RX,Y]+[X,RY]=[X,Y]_R+R[X,Y]$$
is a Lie one. The general condition that the bracket $[\cdot,\cdot]'_R$ is
a Lie one has the form
$$[R^2[X,Y],Z]+c.p.=0.$$
The bracket $[\cdot,\cdot]'_R$ plays a crucial role in the factorization
method [1]. Let us consider an example when $R^2\ne1$, nevertheless,
$[R^2[X,Y],Z]+c.p.=0$.

\remark{Example 1} Let us consider the Witt algebra with the basis
$e_k$ ($k\in\Bbb Z$) and the commutation relations
$$[e_i,e_j]=(i-j)e_{i+j}.$$
Introduce the operators $R_n$ of the form
$$R_n(e_i)=e_{i+n}.$$
Each operator $R_n$ defines a structure of $\mYB$-algebra on the Witt algebra,
for none $n$ the squares of operators $R_n$ are not equal to the identity but
for all $n$ the equality $[R^2_n[X,Y],Z]+c.p.=0$ holds.
\endremark

\remark{Remark 3} The linear $\Gamma$-bunches of Lie algebras do not form
a submanifold (in sense of [4,5]) of the manifold of all linear bunches
of Lie algebras. The identities, which connect the brackets $[\cdot,\cdot]$
and $[\cdot,\cdot]_R$ (as well as $[\cdot,\cdot]$ and $[\cdot,\cdot]'_R$)
are not known.
\endremark

\definition{Definition 1B} A linear bunch of Lie algebras $\frak g$ is called
the {\it linear $\PG$-bunch} if it admits a monomorphism into a linear
$\Gamma$-bunch of Lie algebras.
\enddefinition

As it was marked above a description of the class of all linear $\PG$-bunches
of Lie algebras by the complete system of identities is not known.

Let us consider an example of the linear $\PG$-bunch of Lie algebras,
whcih is not a linear $\Gamma$-bunch.

\remark{Example 2\/ {\rm (cf.[3])}} Let $\frak g=\soa(n)$ be the Lie algebra
of all skew-symmetric $n\times n$ matrices with the standard commutator
$[X,Y]=XY-YX$ and $Q$ be an arbitrary nonscalar symmetric $n\times n$ matrix.
Define a family of compatible Lie brackets $[\cdot,\cdot]_\lambda$ as
$$[X,Y]_\lambda=[X,Y]+\lambda(XQY-YQX).$$
The obtained linear bunch of Lie algebras is not a linear $\Gamma$-bunch,
because in view of the nonscalarity of the matrix $Q$ some of the algebras
$\frak g_\lambda$ being isomorphic to $\soa(p,q)$ do not admit any
homomorphisms into $\soa(n)$. On the other hand this bunch may be imbed
into the bunch of all $n\times n$ matrices with the same brackets, which
is a linear $\Gamma$-bunch with homomorphisms $R_\lambda$ defined by the
operators of the left or right multiplication on the matrices
$Q_\lambda=1+\lambda Q$.
\endremark

\subhead 1.3. Lie $\mYB$-algebras and linear $\Gamma$-bunches of Lie
algebras\endsubhead

\definition{Definition 2A} The {\it Lie $\mYB$-algebra\/} is a Lie
algebra $\frak g$ with the bracket $[\cdot,\cdot]$, supplied by an
operator $R:\frak g\mapsto\frak g$ such that
$$R[RX,Y]+R[X,RY]=[RX,RY]+R^2[X,Y].$$
\enddefinition

\proclaim{Proposition 2} The Jacobi identity holds for the bracket
$$[X,Y]_R=[RX,Y]+[X,RY]-R[X,Y]$$
in any Lie $\mYB$-algebra.
\endproclaim

Theorem 1A states that any linear $\Gamma$-bunch of Lie algebras is a Lie
$\mYB$-algebra. The inverse statement is also true.

\proclaim{Theorem 1B} Any Lie $\mYB$-algebra $(\frak g,R)$ is supplied
by a structure of a linear $\Gamma$-bunch of Lie algebras with brackets
$[\cdot,\cdot]_\lambda=[\cdot,\cdot]+\lambda[\cdot,\cdot]_R$ and
homomorphisms $R_\lambda=1+\lambda R$. In particular, the brackets
$[\cdot,\cdot]$ and $[\cdot,\cdot]_R$ are compatible.
\endproclaim

\proclaim{Proposition 3} For any Lie $\mYB$-algebra $(\frak g,R)$
an operator $R'$ of the form
$f(R)\!=\!a_0\!+\!a_1R\!+\!a_2R^2\!+\!\ldots\!+\!a_nR^n$
defines a structure of the Lie $\mYB$-algebra $(\frak g,f(R))$.
\endproclaim

Hence, the brackets $[\cdot,\cdot]_{R^n}$ ($n\in\Bbb Z_+$) are compatible
between each other. Pro\-po\-sition 3 allows to construct the operators $R'$
such that $(R')^2=1$ for the subsequent applications of the factorization
method.

\subhead 1.4. Lie bi-$\mYB$-algebras\endsubhead

\proclaim{Proposition 4} Let $\frak A$ be an associative algebra then its
commutator algebra $\frak A_{[\cdot,\cdot]}$ is a Lie $\mYB$-algebra, where
$R(X)=R_Q^l(X)=QX$ or $R(X)=R_Q^r(X)=XQ$ are the operators of the left or
right multiplication on the elements $Q$ of the associative algebra $\frak A$.
Moreover,
$$[X,Y]_{R^r_Q}=[X,Y]_{R^l_Q}=XQY-YQX.$$
\endproclaim

This construction motivates the following definition.

\definition{Definition 2B} The {\it Lie bi-$\mYB$-algebra\/} is a Lie algebra
$\frak g$ with bracket $[\cdot,\cdot]$ supplied by two commuting operators
$R_1$ and $R_2$ such that $(\frak g,R_1)$ and $(\frak g,R_2)$ are the Lie
$\mYB$-algebras with identical brackets $[\cdot,\cdot]_{R_1}$ and
$[\cdot,\cdot]_{R_2}$.
\enddefinition

Thus proposition 4 states that the fixing of any element $Q$ of the
associative algebra $\frak A$ supplies this algebra by a structure of the
Lie bi-$\mYB$-algebra $(\frak A_{[\cdot,\cdot]},R^r_Q,R^l_Q)$.

\remark{Remark 4} In any Lie bi-$\mYB$-algebra $(\frak g,R_1,R_2)$ the
operator $R_1\!-\!R_2$ is a derivative of the bracket $[\cdot,\cdot]$ of
the Lie algebra $\frak g$.
\endremark

\remark{Remark 5} A Lie $\mYB$-algebra $(\frak g,R)$ is the Lie
bi-$\mYB$-algebra if and only if there exists a derivative $\xi$ of the Lie
algebra $\frak g$ commuting with $R$ such that
$$[\xi X,\xi Y]=[SX,Y]+[X,SY]-S[X,Y],\qquad S=R\xi,$$
and $R_1=R$, $R_2=R+\xi$.
\endremark

Hence, in the Lie bi-$\mYB$-algebra the operator $\xi$ is a derivative of both
brackets $[\cdot,\cdot]$ and $[\cdot,\cdot]_R$.

\proclaim{Proposition 5} For any Lie bi-$\mYB$-algebra $(\frak g,R_1,R_2)$
and any polynomial $f(x)$ the triple $(\frak g,f(R_1),f(R_2))$ is a Lie
bi-$\mYB$-algebra.
\endproclaim

\remark{Remark 6} In any Lie bi-$\mYB$-algebra $(\frak A_{[\cdot,\cdot]},
R^r_Q,R^l_Q)$ the bracket $[\cdot,\cdot]^q$ defined as
$$[X,Y]^q=[R^r_QX,R^l_QY]+[R^l_QX,R^r_QY]-R^r_QR^l_Q[X,Y],$$
obeys the Jacobi identity
$$[[X,Y]^q,Z]^q+c.p.=0,$$
and also is compatible with $[\cdot,\cdot]$, $[\cdot,\cdot]_{R_1}$ and
$[\cdot,\cdot]_{R_2}$. Moreover, the identity
$$[X,Y]^q=[X,Y]_{R_1^2}=[X,Y]_{R_2^2}$$
holds for all $X$ and $Y$.
\endremark

This fact motivates the following definition.

\definition{Definition 2C} A Lie bi-$\mYB$-algebra $(\frak g,R_1,R_2)$
is called {\it even-tempered} if the identities
$$\aligned
[R_1X,R_2Y]+[R_2X,R_1Y]-R_1R_2[X,Y]&=[R_1^2X,Y]+[X,R_1^2,Y]-R_1^2[X,Y],\\
[R_1X,R_2Y]+[R_2X,R_1Y]-R_1R_2[X,Y]&=[R_2^2X,Y]+[X,R_2^2,Y]-R_2^2[X,Y]
\endaligned
$$
hold.
\enddefinition

Note that one may claim only one of identities in definition 2C in view of
proposition 5.

\remark{Remark 7} The identities in the even-tempered algebra have the form
$$[RX,\xi Y]+[\xi X,RY]-R\xi[X,Y]=[R^2X,Y]-2[RX,RY]+[X,R^2Y]$$
in terms of $R$ and $\xi$.
\endremark

Thus, remark 6 means that the Lie bi-$\mYB$-algebra $(\frak
A_{[\cdot,\cdot]},R^r_Q,R^l_Q)$ constructed from an associative algebra
$\frak A$ and its element $Q$ in proposition 4 is even-tempered.

\head\bf 2. Representations of linear bunches of Lie algebras
\endhead

This paragraph is devoted to representations of linear bunches of Lie
algebras, which naturally appear in the quantization of linearly depending
on a parameter and compatible (linear) Poisson brackets (such situation is
usual for systems in external magnetic fields [6]).

\definition{Definition 3 {\rm (cf.[3])}} The {\it representation of the
linear $\PG$-bunch $\frak g$ of Lie algebras in the linear space $H$}
is a mapping $T:\frak g\mapsto\End(H)$ for which there exists an operator
$Q_R$ in $H$ such that
$$T([X,Y]_\lambda)=T(X)(1+\lambda Q_R)T(Y)-T(Y)(1+\lambda Q_R)T(X)$$
for all $X$ and $Y$ from $\frak g$. A linear $\PG$-bunch of Lie algebras
is called {\it representable} if it admits a faithful representation.
\enddefinition

The claim of a representability of a linear $\PG$-bunch is rather strong.
Let us investigate some necessary conditions of representability
of linear $\PG$-bunches.

Define for two Lie brackets $[\cdot,\cdot]_\alpha$ and
$[\cdot,\cdot]_\beta$ in the linear space $V$ their $\lozenge$-product
(depending on the element of $V$) as
$$\aligned
[X,Y]_{\alpha\underset Z\to\lozenge\beta}=&\frac12\left(
[[X,Z]_\alpha,Y]_\beta+[[X,Y]_\alpha,Z]_\beta+[[Z,Y]_\alpha,X]_\beta-\right.\\
&\left.-[[X,Z]_\beta,Y]_\alpha-[[X,Y]_\beta,Z]_\alpha-[[Z,Y]_\beta,X]_\alpha
\right).\endaligned
$$
Let us formulate the main theorem on the necessary conditions of the
rep\-re\-sen\-ta\-bi\-li\-ty of linear $\Gamma$-bunches of Lie algebras.

\proclaim{Theorem 2} The linear $\PG$-bunch of Lie algebras $\frak g$ is
representable only if there exists a linear family of (compatible) Lie
brackets, which contains the prescribed bunch and which is closed under
the $\lozenge$-product. If such family $V$ of Lie brackets exists then
$(\frak g,V)$ possesses a natural structure of an isotopic pair.
\endproclaim

One has the following corollary for linear $\PG$-bunches of Lie algebras.

\proclaim{Corollary} A linear $\Gamma$-bunch of Lie algebras defined by the
Lie $\mYB$-algebra $(\frak g,R)$ is representable only if the pair
$(\frak g,R+\lambda[\ad Z,R])$ is a Lie $\mYB$-algebra for any element
$Z$ of the Lie algebra $\frak g$ and an arbitrary number $\lambda$.
\endproclaim

Note that each homomorphism of the Lie $\mYB$-algebra $(\frak g,R)$ into
the Lie bi-$\mYB$-algebra $(\frak A_{[\cdot,\cdot]},R^r_Q,R^l_Q)$ of
proposition 4 allows to construct representations of the linear
$\Gamma$-bunch $\frak g$ in the spaces of representations of the associative
algebra $\frak g$. The inverse statement is not true in general, i.e.
a representation $T$ of the linear $\Gamma$-bunch $\frak g$ in the linear
space $H$ does not obligatory define a homomorphism of the Lie $\mYB$-algebra
$(\frak g,R)$ into the Lie bi-$\mYB$-algebra $(\Mat_{[\cdot,\cdot]}(n),R^r_Q,
R^l_Q)$ (i.e. a homomorphism of $(\frak g,R)$ into
$(\Mat_{[\cdot,\cdot]}(n),R^l_Q)$ as algebras with operators) that is
demonstarted by the following example.

\remark{Example 3} Let us consider the Lie algebra $\sla(2,\Bbb C)$ with the
basis $L_{-1}$, $L_0$, $L_1$ and commutation relations
$[L_i,L_j]=(i-j)L_{i+j}$ and the operator $R$ on it: $RL_i=iL_i$. The pair
$(\sla(2,\Bbb C),R)$ is a Lie $\mYB$-algebra and defines the linear
$\Gamma$-bunch $\frak g$. The two-dimensional fundamental representation of
the Lie algebra $\sla(2,\Bbb C)$ realizes a representation of this bunch
with $Q=L_0$, however, none homomorphism of the Lie $\mYB$-algebra
$(\sla(2,\Bbb C),R)$ into the Lie bi-$\mYB$-algebra $(\Mat_{[\cdot,\cdot]}(n),
R^r_Q,R^l_Q)$ corresponds to it, because the equality $T(R)T(L_0)=0$ implies
$T(R)=0$ in view of the invertibility of the operator $T(L_0)$.
\endremark

\Refs
\roster
\item"[1]" {\it Semenov-Tyan-Shanskii M.A.}, What the classical $r$-matrix
is. Funkts.anal.i ego prilozh. 17(4) (1983) 17-33.
\item"[2]" {\it Gelfand I.M., Dorfman I.Ya.}, Hamiltonian operators and
classical Yang-Baxter equation. Funkts.anal.i ego prilozh. 16(4) (1982) 1-9.
\item"[3]" {\it Juriev D.V.}, Topics in isotopic pairs and their
representations. I,II. Teor.Matem.Fiz. 105(1) (1995) 18-28; 111(1) (1997)
149-158.
\item"[4]" {\it Kon P.}, Universal algebra. Moscow, 1968.
\item"[5]" {\it Mal'tsev A.I.}, Algebraic systems. Moscow, 1970.
\item"[6]" {\it Dubrovin B.A., Krichever I.M., Novikov S.P.}, Integrable
systems. In ``Current Math. Problems. Basic Directions IV''. Moscow, VINITI,
1988.
\endroster
\endRefs
\newpage
\centerline{\bf APPENDIX: THE ORIGINAL RUSSIAN VERSION OF ARTICLE}
\ \newline
\centerline{\cyb IZOTOPIChESKIE PARY I IH PREDSTAVLENIYa. {\bf III}.}
\centerline{\cyb PUChKI ALGEBR LI I MODIFITsIROVANNOE}
\centerline{\cyb URAVNENIE YaNGA-BAKSTERA}
\ \newline
\centerline{\cyb D.V.Yurp1ev}
\ \newline\eightpoint
\centerline{\cyre Tsentr matematicheskoe0 fiziki i informatiki ``Talassa
E1teriya'',}
\centerline{\cyre ul.Miklukho-Maklaya 20-180, Moskva 117437 Rossiya.}
\centerline{E-mail: denis\@juriev.msk.ru}
\ \tenpoint\newline
\centerline{q-alg/9708027}
\ \newline

\cyr Dannaya nebolp1shaya zametka yavlyaet\-sya kratkim kommentariem k
ne\-ko\-to\-rym aspektam teorii klassicheskoe0 $r$-matritsy [1] i
bigamilp1tonova for\-ma\-liz\-ma [2], motivirovkoe0 kotorogo sluzhat
konstruktsii rabot [3].

\head\cyb 1. Linee0nye puchki algebr Li i $\mYB$-algebry Li
\endhead

\cyr Dannye0 paragraf posvyashchen analizu nekotoryh algebraicheskih
struk\-tur, svyazannyh s linee0nymi puchkami algebr Li.

\subhead\cyb 1.1. $\Gamma$-puchki algebr Li i kasatelp1naya skobka\endsubhead

\definition{\cyb Opredelenie 1A}\cyr {\cyi Puchkom (odnoparametricheskim)
algebr Li\/} nazyvaet\-sya linee0noe prostranstvo $\frak g$, snabzhennoe
odnoparametricheskim semee0stvom lievskih skobok $[\cdot,\cdot]_\lambda$
($\lambda\in\Bbb R$). Algebra Li iz puchka $\frak g$, opredelyaemaya
skob\-koe0 $[\cdot,\cdot]_\lambda$, oboznachaet\-sya $\frak g_\lambda$.
Puchok algebr Li $\frak g$ nazyvaet\-sya {\cyi $\Gamma$-puchkom}, esli
zadany gomomorfizmy $R_\lambda$ iz algebr Li $\frak g_\lambda$ v algebru Li
$\frak g_0$, inymi slovami, dlya vseh $X$ i $Y$ iz $\frak g$ vypolneno
ravenstvo $$R_\lambda[X,Y]_\lambda=[R_\lambda X,R_\lambda Y],$$
gde $[\cdot,\cdot]=[\cdot,\cdot]_0$.
\enddefinition

\cyr V dalp1nee0shem budut rassmatrivatp1sya gladkie puchki, t.e. puchki
dlya kotoryh lievskie skobki obrazuyut gladkoe semee0stvo, i gladkie
$\Gamma$-puchki, dlya kotoryh gladkim yavlyaet\-sya takzhe semee0stvo
gomomorfizmov $R_\lambda$.

\proclaim{\cyb Predlozhenie 1}\cyi Pustp1 $\frak g$ - gladkie0
$\Gamma$-puchok algebr Li, opredelim ka\-sa\-telp1\-nuyu skobku
$[\cdot,\cdot]_R$ sleduyushchim obrazom:
$$[X,Y]_R=\left.\frac{d[X,Y]_\lambda}{d\lambda}\right|_{\lambda=0}.$$
Togda
$$[X,Y]_R=[RX,Y]+[X,RY]-R[X,Y],$$
gde
$$R=\left.\frac{dR_\lambda}{d\lambda}\right|_{\lambda=0}.$$
\endproclaim

\remark{\cyi Zamehcanie 1}\cyr Skobka $[\cdot,\cdot]_R$ trivialp1na
(tozhdestvenno ravna nulyu) togda i tolp1ko togda, kogda $R$ --
differentsirovanie skobki $[\cdot,\cdot]$.
\endremark

\remark{\cyi Zamechanie 2}\cyr Uslovie lievosti skobki $[\cdot,\cdot]_R$
imeet vid
$$\bigl([B(X,Y),Z]+B([X,Y],Z)\bigr)+c.p.=0,$$
gde
$$B(X,Y)=R[RX,Y]+R[X,RY]-[RX,RY]-R^2[X,Y].$$
Razvernutoe vyrazhenie dlya $B(X,Y)$ mozhet bytp1 skomponovano v
dvu\-chlen\-nuyu formulu:
$$B(X,Y)=R[X,Y]_R-[RX,RY].$$
\endremark

\subhead\cyb 1.2. Linee0nye $\Gamma$-puchki algebr Li i modifitsirovannoe
uravnenie Yanga-Bakstera\endsubhead\cyr Otmetim, chto skobka $[\cdot,\cdot]$
i kasatelp1naya skobka $[\cdot,\cdot]_R$ v proizvolp1nom $\Gamma$-puchke
algebr Li vsegda udovletvoryaet tozhdestvu
$$\bigl([[X,Y]_R,Z]+[[X,Y],Z]_R\bigr)+c.p.=0.$$
Kak sledstvie, esli kasatelp1naya skobka $[\cdot,\cdot]_R$ lievskaya, to ona
soglasovana so skobkoe0 $[\cdot,\cdot]$. E1to obstoyatelp1stvo yavlyaet\-sya
motivirovkoe0 dlya ras\-smot\-re\-niya linee0nyh $\Gamma$-puchkov algebr Li.

Puchok algebr Li $\frak g$ nazyvaet\-sya linee0nym, esli lievskie skobki
$[\cdot,\cdot]_\lambda$ obrazuyut linee0noe semee0stvo, inymi slovami
$[\cdot,\cdot]_\lambda=(1-\lambda)[\cdot,\cdot]+\lambda[\cdot,\cdot]_1$.
$\Gamma$-puchok algebr Li nazyvaet\-sya linee0nym, esli ne tolp1ko lievskie
skobki, no i gomomorfizmy $R_\lambda$ obrazuyut linee0noe semee0stvo, t.e.
$R_\lambda=1+\lambda R$. Trebovanie linee0nosti semee0stva gomomorfizmov
$R_\lambda$ estestvenno, poskolp1ku dlya nelinee0nogo semee0stva
$R_\lambda$ ego linee0naya chastp1 $1+\lambda R$ ($R=\left.\frac{dR_\lambda}
{d\lambda}\right|_{\lambda=0}$) zadaet to zhe linee0noe semee0stvo skobok
$[\cdot,\cdot]_\lambda$.

\proclaim{\cyb Teorema 1A}\cyi Dlya linee0nogo $\Gamma$-puchka algebr Li
$\frak g$ kasatelp1naya skobka $[\cdot,\cdot]_R$ (kotoraya v e1tom sluchae
lievskaya) udovletvoryaet usloviyu
$$R[X,Y]_R=[RX,RY],$$
inymi slovami,
$$R[RX,Y]+R[X,RY]=[RX,RY]+R^2[X,Y].$$
\endproclaim

\cyr Teorema oznachaet, chto uslovie lievosti skobki $[\cdot,\cdot]_R$ iz
zamechaniya 2 v sluchae linee0nyh $\Gamma$-puchkov faktorizuet\-sya i
svodit\-sya k trebovaniyu tozhdestvennogo ravenstva vyrazheniya $B(X,Y)$
nulyu. Otmetim, chto li\-e\-vostp1 skobki $[\cdot,\cdot]_R$ v silu ee
soglasovannosti so skobkoe0 $[\cdot,\cdot]$ oznachaet, chto puchok
$[\cdot,\cdot]+\lambda[\cdot,\cdot]_R$ -- linee0nye0 puchok algebr Li,
odnako, on ne obyazan bytp1 $\Gamma$-puchkom, poe1tomu v obshchem sluchae
uslovie lievosti iz zamechaniya 2 ne faktorizuet\-sya.

Poslednee tozhdestvo v teoreme 1A sovpadaet s modifitsirovannym
klas\-si\-ches\-kim uravneniem Yanga-Bakstera [1] pri $R^2=1$, chto
vypolnyaet\-sya dlya naibolee vazhnyh primerov, rassmotrennyh v [1].
V e1tom sluchae, kak pokazano v [1], lievskoe0 skobkoe0 takzhe
yavlyaet\-sya skobka $[\cdot,\cdot]'_R$, zadavaemaya sleduyushchim obrazom
$$[X,Y]'_R=[RX,Y]+[X,RY]=[X,Y]_R+R[X,Y].$$
Obshchee uslovie lievosti skobki $[\cdot,\cdot]'_R$ vyglyadit sleduyushchim
obrazom:
$$[R^2[X,Y],Z]+c.p.=0.$$
Skobka $[\cdot,\cdot]'_R$ igraet klyuchevuyu rlp1 v metode faktorizatsii [1].
Privedem primer, kogda $R^2\ne1$, odnako, $[R^2[X,Y],Z]+c.p.=0$.

\remark{\cyi Primer 1}\cyr Rassmotrim algebru Vitta s bazisom $e_k$
($k\in\Bbb Z$) i kom\-mu\-ta\-tsi\-on\-ny\-mi sootnosheniyami
$$[e_i,e_j]=(i-j)e_{i+j}.$$
Vvedem operatory $R_n$ sleduyushchego vida
$$R_n(e_i)=e_{i+n}.$$
Kazhdye0 iz operatorov $R_n$ zadaet na algebre Vitta strukturu
$\mYB$-al\-geb\-ry, pri e1tom ni dlya kakogo $n$ kvadraty operatorov $R_n$ ne
ravny edinitse, odnako, pri vseh $n$ tozhdestvo $[R^2_n[X,Y],Z]+c.p.=0$
vypolneno.
\endremark

\remark{\cyi Zamechanie 3}\cyr Linee0nye $\Gamma$-puchki algebr Li ne
obrazuyut podmnogoobraziya (v smysle [4,5]) v mnogoobrazii vseh lineenyh
puchkov algebr Li. Tozh\-dest\-va, kotorye svyazyvayut skobki $[\cdot,\cdot]$
i $[\cdot,\cdot]_R$ (kak, vprochem, i $[\cdot,\cdot]$ i $[\cdot,\cdot]'_R$)
neizvestny.
\endremark

\definition{\cyb Opredelenie 1B}\cyr Linee0nye0 puchok algebr Li $\frak g$
nazyvaet\-sya {\cyi linee0nym $\PG$-puchkom}, esli on dopuskaet
monomorfizm v nekotorye0 linee0nye0 $\Gamma$-puchok algebr Li.
\enddefinition

\cyr Kak otmechalosp1 vyshe, opisanie klassa linee0nyh $\PG$-puchkov algebr
Li s pomoshchp1yu tozhdestv neizvestno.

Privedem primer linee0nogo $\PG$-puchka algebr Li, ne yavlyayushchegosya
li\-nee0\-nym $\Gamma$-puchkom.

\remark{\cyi Primer 2\/ {\cyr (sr.[3])}}\cyr Pustp1 $\frak g=\soa(n)$ --
algebra Li vseh kososimmetricheskih matrits $n\times n$ so standartnym
kommutatorom $[X,Y]=XY-YX$, i $Q$ -- proizvolp1naya neskalyarnaya
simmetricheskaya matritsa $n\times n$. Zadadim semee0stvo soglasovannyh
lievskih skobok $[\cdot,\cdot]_\lambda$ sleduyushchim obrazom
$$[X,Y]_\lambda=[X,Y]+\lambda(XQY-YQX).$$
Poluchennye0 linee0nye0 puchok algebr Li ne yavlyaet\-sya linee0nym
$\Gamma$-puchkom, poskolp1ku v silu neskalyarnosti matritsy $Q$ nekotorye
iz algebr $\frak g_\lambda$, bu\-du\-chi izomorfnymi $\soa(p,q)$, ne dopuskayut
gomomorfizmov v $\soa(n)$. S drugoe0 storony, dannye0 puchok vkladyvaet\-sya
v puchok vseh matrits $n\times n$ s temi zhe skobkami, kotorye0 yavlyaet\-sya
linee0nym $\Gamma$-puchkom s gomomorfizmami $R_\lambda$, zadavaemymi
operatorami umnozheniya sprava ili sleva na matritsy $Q_\lambda=1+\lambda Q$.
\endremark

\subhead\cyb 1.3. $\mYB$-algebry Li i linee0nye $\Gamma$-puchki algebr
Li\endsubhead

\definition{\cyb Opredelenie 2A}\cyr {\cyi $\mYB$-algebroe0 Li\/}
nazyvaet\-sya algebra Li $\frak g$ so skobkoe0 $[\cdot,\cdot]$,
snabzhennaya operatorom $R:\frak g\mapsto\frak g$ takim, chto
$$R[RX,Y]+R[X,RY]=[RX,RY]+R^2[X,Y].$$
\enddefinition

\proclaim{\cyb Predlozhenie 2}\cyi V $\mYB$-algebre Li dlya skobki
$$[X,Y]_R=[RX,Y]+[X,RY]-R[X,Y]$$
vypolnyaet\-sya tozhdestvo Yakobi.
\endproclaim

\cyr Teorema 1A utverzhdaet, chto vsyakie0 linee0nye0 $\Gamma$-puchok algebr
Li yav\-lya\-et\-sya $\mYB$-algebroe0 Li. Imeet mesto i obratnoe utverzhdenie.

\proclaim{\cyb Teorema 1B}\cyi Vsyakaya $\mYB$-algebra Li $(\frak g,R)$
snabzhaet\-sya strukturoe0 li\-nee0\-no\-go $\Gamma$-puchka algebr Li so
skobkami $[\cdot,\cdot]_\lambda=[\cdot,\cdot]+\lambda[\cdot,\cdot]_R$ i
gomomorfizmami $R_\lambda=1+\lambda R$. V chastnosti, skobki $[\cdot,\cdot]$
i $[\cdot,\cdot]_R$ soglasovany.
\endproclaim

\proclaim{\cyb Predlozhenie 3}\cyi Dlya proizvolp1noe0 $\mYB$-algebry Li
$(\frak g,R)$ lyuboe0 operator $R'$ vida
$f(R)\!=\!a_0\!+\!a_1R\!+\!a_2R^2\!+\!\ldots\!+\!a_nR^n$ zadaet strukturu
$\mYB$-algebry Li $(\frak g,f(R))$.
\endproclaim

\cyr Kak sledstvie, skobki $[\cdot,\cdot]_{R^n}$ ($n\in\Bbb Z_+$) soglasovany
mezhdu soboe0. Pred\-lo\-zhe\-nie 3 pozvolyaet stroitp1 operatory $R'$ takie,
chto $(R')^2=1$, dlya posleduyushchego primeneniya metoda faktorizatsii.

\subhead\cyb 1.4. Bi-$\mYB$-algebry Li\endsubhead

\proclaim{\cyb Predlozhenie 4}\cyi Pustp1 $\frak A$ -- assotsiativnaya
algebra, togda kom\-mu\-ta\-tor\-naya algebra $\frak A_{[\cdot,\cdot]}$
yavlyaet\-sya $\mYB$-algebroe0 Li, gde $R(X)=R_Q^l(X)=QX$ ili
$R(X)=R_Q^r(X)=XQ$ -- operatory umnozheniya sleva i sprava na e1le\-men\-ty
$Q$ assotsiativnoe0 algebry $\frak A$. Pri e1tom,
$$[X,Y]_{R^r_Q}=[X,Y]_{R^l_Q}=XQY-YQX.$$
\endproclaim

\cyr E1ta konstruktsiya yavlyaet\-sya motivirovkoe0 dlya sleduyushchego
opredeleniya.

\definition{\cyb Opredelenie 2B}\cyr {\cyi Bi-$\mYB$-algebroe0 Li\/}
nazyvaet\-sya algebra Li $\frak g$ so skob\-koe0 $[\cdot,\cdot]$, snabzhennaya
dvumya kommutiruyushchimi operatorami $R_1$ i $R_2$ ta\-ki\-mi, chto
$(\frak g,R_1)$ i $(\frak g,R_2)$ -- $\mYB$-algebry Li s sovpadayushchimi
skobkami $[\cdot,\cdot]_{R_1}$ i $[\cdot,\cdot]_{R_2}$.
\enddefinition

\cyr Takim obrazom, predlozhenie 4 utverzhdaet, chto fiksirovanie
pro\-iz\-volp1\-no\-go e1lementa $Q$ assotsiativnoe0 algebry $\frak A$
snabzhaet ukazannuyu algebru strukturoe0 bi-$\mYB$-algebry Li
$(\frak A_{[\cdot,\cdot]},R^r_Q,R^l_Q)$.

\remark{\cyi Zamechanie 4}\cyr Operator $R_1\!-\!R_2$ v bi-$\mYB$-algebre
$(\frak g,R_1,R_2)$ yavlyaet\-sya dif\-fe\-ren\-tsi\-ro\-va\-ni\-em skobki
$[\cdot,\cdot]$ algebry Li $\frak g$.
\endremark

\remark{\cyi Zamechanie 5}\cyr $\mYB$-algebra Li $(\frak g,R)$ yavlyaet\-sya
bi-$\mYB$-algebroe0 togda i tolp1ko togda, kogda sushchestvuet
differentsirovanie $\xi$ algebry Li $\frak g$, kom\-mu\-ti\-ru\-yu\-shchee
s $R$, takoe, chto
$$[\xi X,\xi Y]=[SX,Y]+[X,SY]-S[X,Y],\qquad S=R\xi.$$
Pri e1tom $R_1=R$, $R_2=R+\xi$.
\endremark

\cyr Kak sledstvie, operator $\xi$ v bi-$\mYB$-algebre Li yavlyaet\-sya
dif\-fe\-ren\-tsi\-ro\-va\-ni\-em obeih skobok $[\cdot,\cdot]$ i
$[\cdot,\cdot]_R$.

\proclaim{\cyb Predlozhenie 5}\cyi Dlya lyuboe0 bi-$\mYB$-algebry Li
$(\frak g,R_1,R_2)$ i mnogochlena $f(x)$ troe0ka $(\frak g,f(R_1),f(R_2))$
yavlyaet\-sya bi-$\mYB$-algebroe0 Li.
\endproclaim

\remark{\cyi Zamechanie 6}\cyr V bi-$\mYB$-algebre Li
$(\frak A_{[\cdot,\cdot]},R^r_Q,R^l_Q)$ skobka $[\cdot,\cdot]^q$,
zadavaemaya sleduyushchim obrazom
$$[X,Y]^q=[R^r_QX,R^l_QY]+[R^l_QX,R^r_QY]-R^r_QR^l_Q[X,Y],$$
udovletvoryaet tozhdestvu Yakobi
$$[[X,Y]^q,Z]^q+c.p.=0,$$
a takzhe soglasovana so skobkami $[\cdot,\cdot]$, $[\cdot,\cdot]_{R_1}$ i
$[\cdot,\cdot]_{R_2}$. Bolee togo, dlya vseh $X$ i $Y$ vypolneno ravenstvo
$$[X,Y]^q=[X,Y]_{R_1^2}=[X,Y]_{R_2^2}.$$
\endremark

\cyr E1tot fakt yavlyaet\-sya motivirovkoe0 dlya sleduyushchego opredeleniya.

\definition{\cyb Opredlenie 2V}\cyr Bi-$\mYB$-algebra Li $(\frak g,R_1,R_2)$
nazyvaet\-sya {\cyi urav\-no\-ve\-shen\-noe0}, esli vypolneny sleduyushchie
tozhdestva
$$\aligned
[R_1X,R_2Y]+[R_2X,R_1Y]-R_1R_2[X,Y]&=[R_1^2X,Y]+[X,R_1^2,Y]-R_1^2[X,Y],\\
[R_1X,R_2Y]+[R_2X,R_1Y]-R_1R_2[X,Y]&=[R_2^2X,Y]+[X,R_2^2,Y]-R_2^2[X,Y].
\endaligned
$$
\enddefinition

\cyr Otmetim, chto v opredelenii 2V dostatochno potrebovatp1 vypolneniya
odnogo iz tozhdestv v silu predlozheniya 5.

\remark{\cyi Zamechanie 7}\cyr Tozhdestva uravnoveshennosti v terminah $R$
i $\xi$ imeyut vid
$$[RX,\xi Y]+[\xi X,RY]-R\xi[X,Y]=[R^2X,Y]-2[RX,RY]+[X,R^2Y].$$
\endremark

\cyr Takim obrazom, zamechanie 6 oznachaet, chto bi-$\mYB$-algebra $(\frak
A_{[\cdot,\cdot]},R^r_Q,R^l_Q)$, postroennaya v predlozhenii 4 po
assotsiativnoe0 algebre $\frak A$ i ee e1lementu $Q$, yavlyaet\-sya
uravnoveshennoe0.

\head\cyb 2. Predstavleniya linee0nyh puchkov algebr Li
\endhead

\cyr Dannye0 paragraf posvyashchen predstavleniyam linee0nyh puchkov
algebr Li, kotorye estestvenno voznikayut v ramkah kvantovaniya linee0no
za\-vi\-sya\-shchih ot parametra soglasovannyh (linee0nyh) skobok Puassona
(chto harakterno dlya sistem vo vneshnem magnitnom pole [6]).

\definition{\cyb Opredelenie 3 {\cyr (sr.[3])}}\cyr {\cyi Predstavleniem
linee0nogo $\PG$-puchka $\frak g$ algebr Li v linee0nom prostranstve $H$},
nazyvaet\-sya otobrazhenie $T:\frak g\mapsto\End(H)$, dlya kotorogo
sushchestvuet operator $Q_R$ v $H$ takoe0, chto
$$T([X,Y]_\lambda)=T(X)(1+\lambda Q_R)T(Y)-T(Y)(1+\lambda Q_R)T(X),$$
dlya proizvolp1nyh $X$ i $Y$ iz $\frak g$. Linee0nye0 $\PG$-puchok algebr Li
nazyvaet\-sya {\cyi predstavimym}, esli on obladaet hotya by odnim tochnym
predstavleniem.
\enddefinition

\cyr Trebovanie predstavimosti linee0nogo $\PG$-puchka dostatochno
silp1noe. Izuchim ryad neobhodimyh uslovie0 predstavimosti linee0nyh
$\PG$-puchkov.

Opredelim dlya dvuh lievskih skobok $[\cdot,\cdot]_\alpha$ i
$[\cdot,\cdot]_\beta$, zadannyh v linee0nom prostranstve $V$, ih
$\lozenge$-proizvedenie (zavisyashchee ot e1lementa prostranstva $V$)
sleduyushchim obrazom:
$$\aligned
[X,Y]_{\alpha\underset Z\to\lozenge\beta}=&\frac12\left(
[[X,Z]_\alpha,Y]_\beta+[[X,Y]_\alpha,Z]_\beta+[[Z,Y]_\alpha,X]_\beta-\right.\\
&\left.-[[X,Z]_\beta,Y]_\alpha-[[X,Y]_\beta,Z]_\alpha-[[Z,Y]_\beta,X]_\alpha
\right).\endaligned
$$
Sformuliruem osnovnuyu teoremu o neobhodimyh usloviyah pred\-sta\-vi\-mos\-ti
linee0nyh $\Gamma$-puchkov algebr Li.

\proclaim{\cyb Teorema 2}\cyi Linee0nye0 $\PG$-puchok algebr Li $\frak g$
predstavim, tolp1ko esli su\-shchest\-vu\-et linee0noe semee0stvo
(sovmestimyh) lievskih skobok, so\-der\-zha\-shchee ukazannye0 puchok i
zamknutoe otnositelp1no $\lozenge$-proizvedeniya. Esli ta\-koe semee0stvo $V$
lievskih skobok sushchestvuet, to $(\frak g,V)$ nadelyaet\-sya
es\-test\-ven\-noe0 strukturoe0 izotopicheskoe0 pary.
\endproclaim

\cyr Dlya linee0nyh $\PG$-puchkov algebr Li imeem sleduyushchee sledstvie.

\proclaim{\cyb Sledstvie}\cyi Linee0nye0 $\Gamma$-puchok algebr Li,
zadavaemye0 $\mYB$-algebroe0 Li\linebreak $(\frak g,R)$, predstavim, tolp1ko
esli dlya lyubogo e1lementa $Z$ algebry Li $\frak g$ i pro\-iz\-volp1\-no\-go
chisla $\lambda$ para $(\frak g,R+\lambda[\ad Z,R])$ yavlyaet\-sya
$\mYB$-algebroe0 Li.
\endproclaim

\cyr Otmetim, chto po vsyakomu gomomorfizmu $\mYB$-algebry Li $(\frak g,R)$
v bi-$\mYB$-algebru Li $(\frak A_{[\cdot,\cdot]},R^r_Q,R^l_Q)$ predlozheniya
4 stroyat\-sya predstavleniya linee0nogo $\Gamma$-puchka $\frak g$ v
prostranstvah predstavlenie0 assotsiativnoe0 al\-geb\-ry $\frak g$. Obratnoe,
voobshche govorya, ne verno, t.e. predstavlenie $T$ linee0nogo
$\Gamma$-puchka $\frak g$ v linee0nom prostranstve $H$ ne vsegda
opredelyaet gomomorfizm $\mYB$-algebry Li $(\frak g,R)$ v bi-$\mYB$-algebru
Li $(\Mat_{[\cdot,\cdot]}(n),R^r_Q,R^l_Q)$ (t.e. go\-mo\-mor\-fizm
$(\frak g,R)$ v $(\Mat_{[\cdot,\cdot]}(n),R^l_Q)$ kak algebr s operatorami),
chto po\-ka\-zy\-va\-et sleduyushchie0 primer.

\remark{\cyi Primer 3}\cyr Rassmotrim algebru Li $\sla(2,\Bbb C)$ s bazisom
$L_{-1}$, $L_0$, $L_1$ i kom\-mu\-ta\-tsi\-on\-ny\-mi sootnosheniyami
$[L_i,L_j]=(i-j)L_{i+j}$, a takzhe operator $R$ na nee0: $RL_i=iL_i$.
Para $(\sla(2,\Bbb C),R)$ yavlyaet\-sya $\mYB$-algebroe0 Li i opredelyaet
linee0nye0 $\Gamma$-puchok $\frak g$. Fundamentalp1noe dvumernoe
predstavlenie algebry Li $\sla(2,\Bbb C)$ osushchestvlyaet predstavlenie
e1togo puchka s $Q=L_0$, odnako, emu ne sootvetstvuet nikakoe0 gomomorfizm
$\mYB$-algebry Li $(\sla(2,\Bbb C),R)$ v bi-$\mYB$-algebru Li
$(\Mat_{[\cdot,\cdot]}(n),R^r_Q,R^l_Q)$ poskolp1ku v silu obratimosti
operatora $T(L_0)$ ravenstvo $T(R)T(L_0)=0$ vlechet $T(R)=0$.
\endremark

\Refs\nofrills{\cyb Spisok literatury}
\roster
\item"[1]" {\cyie Semenov-Tyan-Shanskie0 M.A.}, {\cyre Chto takoe
klassicheskaya $r$-matritsa} /\!/ {\cyre Funkts.anal. i ego prilozh. 1983.
T.17(4). S.17-33.}
\item"[2]" {\cyie Gelp1fand I.M., Dorfman I.Ya.}, {\cyre Gamilp1tonovy
operatory i klassicheskoe urav\-ne\-nie Yanga-Bakstera} /\!/
{\cyre Funkts.anal.i ego prilozh. 1982. T.16(4). S.1-9}.
\item"[3]" {\cyie Yurp1ev D.V.}, {\cyre Izotopicheskie pary i ih
predstavleniya.} I,II /\!/ {\cyre TMF. 1995. T.105(1). S.18-28};
{\cyre TMF. 1997. T.111(1). S.149-158.}
\item"[4]" {\cyie Kon P.}, {\cyre Universalp1naya algebra. M., 1968.}
\item"[5]" {\cyie Malp1tsev A.I.}, {\cyre Algebraicheskie sistemy. M., 1970.}
\item"[6]" {\cyie Dubrovin B.A., Krichever I.M., Novikov S.P.},
{\cyre Integriruemye sistemy / Sov\-re\-men\-nye probl.matematiki. Fundam.
napravleniya {\rm IV}. M., VINITI, 1988.}
\endroster
\endRefs
\enddocument